\documentclass[journal]{IEEEtran}

\usepackage{url}
\usepackage{amsfonts}
\usepackage{amsmath,amssymb,graphicx,epsfig}
\usepackage{enumerate,amsthm,epstopdf}
\usepackage{algpseudocode}
\usepackage{algorithmicx}
\usepackage{algorithm}
\usepackage{comment}
\usepackage{color}
\usepackage{wrapfig,lipsum,booktabs}
\usepackage{verbatim}
\usepackage{balance}

\DeclareMathOperator{\prox}{prox}

\newcommand{\argmin}{\operatornamewithlimits{argmin}}

\def\x{{\mathbf x}}
\def\d{{\mathbf d}}
\def\y{{\mathbf y}}

\def\V{{\mathbf V}}
\def\A{{\mathbf A}}
\def\M{{\mathbf M}}
\def\Q{{\mathbf Q}}
\def\thetab{{\boldsymbol{\theta}}}

\def\bR{{\mathbb R}}

\def\NPDF{{\mathcal N}}

\def\f0{{\mathbf 0}}

\usepackage{mathtools}
\newtheorem{defn}{Definition}

\theoremstyle{definition}

\newtheorem{rem}{Remark}

\usepackage[noadjust]{cite}

\begin{document}
%
\title{A probabilistic incremental \\ proximal gradient method}

\author{\"Omer Deniz Akyildiz, \'Emilie Chouzenoux, V\'ictor Elvira, Joaqu\'in M\'iguez \\
\thanks{\"O. D. Akyildiz is with the Dept. of Computer Science and Dept. of Statistics at University of Warwick and Alan Turing Institute, London, UK. Email: \texttt{omer.akyildiz@warwick.ac.uk}}
\thanks{\'E. Chouzenoux is with the Center for Visual Computing, INRIA Saclay, CentraleSup\'elec, Gif-sur-Yvette, France.}
\thanks{V. Elvira is with IMT Lille Douai \& CRIStAL laboratory (UMR CNRS 9189), Villeneuve d'Ascq, France.}
\thanks{J. M\'iguez is with the Dept. of Signal Theory and Communications, Universidad Carlos III de Madrid, Legan\'es, Spain, 28912.}
\thanks{{\"O}.~D.~A. is funded by the Lloyds Register Foundation programme on Data Centric Engineering through the London Air Quality project and supported by The Alan Turing Institute for Data Science and AI under EPSRC grant EP/N510129/1. J.~M. acknowledges the support of the Spanish \textit{Agencia Estatal de Investigaci\'on} (TEC2015-69868-C2-1-R ADVENTURE) and the Office of Naval Research (N00014-19-1-2226). V.~E. and \'E~.C. acknowledge the support from the \textit{Agence Nationale de la Recherche} of France under PISCES (ANR-17-CE40-0031-01) and MAJIC
(ANR-17-CE40-0004-01) projects.}%
}

\maketitle

\begin{abstract}
In this paper, we propose a probabilistic optimization method, named \textit{probabilistic incremental proximal gradient} (PIPG) method, by developing a probabilistic interpretation of the incremental proximal gradient algorithm. We explicitly model the update rules of the incremental proximal gradient method and develop a systematic approach to propagate the uncertainty of the solution estimate over iterations. The PIPG algorithm takes the form of Bayesian filtering updates for a state-space model constructed by using the cost function. Our framework makes it possible to utilize well-known exact or approximate Bayesian filters, such as Kalman or extended Kalman filters, to solve large-scale regularized optimization problems.
\end{abstract}

%
\IEEEpeerreviewmaketitle

\section{Introduction}
\label{sec:intro}
In this paper, we are interested in optimization problems of the form
\begin{align}\label{eqMainProblem}
\min_{\thetab\in\mathbb{R}^d} f(\thetab) + g(\thetab),
\end{align}
{with $f(\thetab) = \sum_{k=1}^n f_k(\thetab)$ where, for $1 \leq k \leq n$, $f_k: \bR^d \to \bR$, are {nonlinear least squares} functions i.e., for $1 \leq k \leq n$, $f_k = \frac{1}{2}(y_k - h_k(\cdot))^2$, with $y_k \in \bR$ and $h_k: \bR^d \to \bR$ a nonlinear differentiable mapping}. Moreover, $g: \bR^d \to \bR$ is a twice-differentiable regularizer. Because classical optimization schemes may be inefficient to solve \eqref{eqMainProblem} when $n$ is very large, {stochastic or incremental} optimization {methods} have gained a significant momentum. In particular, the stochastic gradient descent (SGD) \cite{RobbinsMonro} has become widely popular to solve such problems. At each SGD iteration, a mini-batch of component functions is randomly selected and a gradient step with respect to this mini-batch is performed. A number of variants of SGD have been since developed (see \cite{bottou2016optimization,Pereyra16} for a review).

The objective function in Eq. \eqref{eqMainProblem} has a sum structure. Therefore, it opens the door for more efficient algorithms than gradient methods, {such as \textit{proximal splitting} methods} \cite{combettes2011proximal,parikh2014proximal}. In particular, the proximal gradient (PG) {method} minimizes a sum of two terms, one being smooth, by alternating gradient steps on the differentiable one and proximal update on the second, thereby exploiting fully the structure of the cost function. Naturally, stochastic extensions of proximal methods have become increasingly popular in the machine learning literature, see, e.g., \cite{bertsekas2011incremental,bertsekas2011incremental2, rosasco2014convergence, atchade2017perturbed, combettes2015stochastic}. The optimization method in consideration in this paper is known as the \textit{incremental proximal gradient (IPG)} algorithm~\cite{bertsekas2011incremental2}, and can be understood as an incremental version of the stochastic proximal gradient method \cite{rosasco2014convergence,atchade2017perturbed}. Similarly to its batch version PG, the IPG method would solve~\eqref{eqMainProblem} by using the gradient of $g$ and the proximal operator of $f_k$ (or vice versa) at each iteration to move within the parameter space. Therefore, the IPG takes advantage of the structure of the cost function while staying computationally efficient for large $n$. 

In this paper, we propose {a} probabilistic IPG (PIPG) method to solve the problem in Eq. \eqref{eqMainProblem}. The PIPG algorithm reads as an {approximate} inference method in a probabilistic state-space model (SSM), tailored to the loss function. To be specific, it takes the form of an extended Kalman filter (EKF) to infer the hidden states of this SSM. This setting yields a probabilistic interpretation which enables the quantification of the uncertainty of the estimates at any time, {extending our previous work \cite{akyildiz2018proximal} which only focused on the case $g=0$}. The posterior covariance matrix {involved in PIPG updates} plays the role of a variable-metric. Thus, another key advantage of PIPG is to provide an adaptive rule for the metric update within the IPG scheme. {Note that the PIPG method is related to the class of \textit{probabilistic numerical methods} (see, e.g., \cite{hennig2015probabilisticnumerics, diaconis1988bayesian, cockayne2017bayesian}), {extending such methods for solving large-scale optimization problems}. {We mention \cite{ollivier2018online} as a related work, that emphasizes the links} between Kalman filtering and the online natural gradient method, which can be viewed as an SGD within a specific variable metric. {In \cite{probabilisticLMS}, connections between LMS and Kalman filters are exploited to propose a new algorithm.} In \cite{vuckovic2018kalman}, the author proposes {some} variance reduction strategies for SGD, relying on {a} Kalman interpretation. In contrast, {in this work we take advantage} of the structure of the cost function itself and {we} focus on the connection between Kalman and proximal methods.}

The paper is organized as follows. In Section~\ref{sec:background}, we briefly {give} background. In Section~\ref{sec:PIPG}, we introduce the new scheme and the update rules in detail. In Section~\ref{sec:experiments}, {we demonstrate the performance of our method, on a ridge regression} and a nonlinear sparse filter identification problem. We conclude with Section~\ref{sec:conclusions}.

\section{Background}\label{sec:background}
Let us start by defining the proximal operator \cite{bauschke2011convex}\footnote{See also \texttt{\url{http://proximity-operator.net/}}}.
\begin{defn} The proximal operator of a convex, proper, lower semi-continuous function $f:\bR^d \to \bR$ within the metric induced by a symmetric, positive definite (SPD) matrix $\V_0~\in~\bR^{d\times d}$ is defined as,
$
\prox_{f,\V_0}(\thetab_0) = \argmin_{\thetab\in\bR^d} f(\thetab) + \frac{1}{2} \|\thetab - \thetab_0\|_{2,\V_0}^2
$
where $\|\thetab\|_{2,\V} :=  {(\thetab^\top \V^{-1} \thetab)^{1/2}}$ is the Mahalanobis distance.
\end{defn}

Let us now present the {Kalman updates} from \cite{akyildiz2018proximal} {which aim at performing Bayesian inference in the case of} the model $p(\thetab) = \NPDF(\thetab;\overline{\thetab}_0,\V_0)$ and $p(y_k|\thetab) = \NPDF(y_k;\x_k^\top \thetab,\gamma^{-1})$, where $\gamma >0$, $\overline{\thetab}_0 \in \mathbb{R}^d$, $\V_0 \in \mathbb{R}^{d \times d}$ SPD, {$\x_k \in \mathbb{R}^d$, for $k = 1,\ldots,n$, are predefined values}, and $(y_k,\thetab)$ are random variables in $\mathbb{R}$ and $\mathbb{R}^d$, respectively. {For this model, assuming that the inputs $\x_{1:k}$ are fixed and the likelihood factorizes as $p(y_{1:k}|\thetab) = \prod_{k=1}^n p(y_k|\thetab)$ (i.e., the observations are conditionally independent)}, the mean and the covariance of the Gaussian posterior $p(\thetab|y_{1:k})=\mathcal{N}(\thetab;\overline{\thetab}_k,\V_k)$ can be written as \cite{akyildiz2018proximal}
\begin{align}
\overline{\thetab}_k &= \overline{\thetab}_{k-1} + \frac{\V_{k-1} \x_k (y_k - \x_k^\top \overline{\thetab}_{k-1})}{\gamma^{-1} + \x_k^\top \V_{k-1} \x_k}, \label{eq:KalmanUpdateMean}\\
\V_k &= \V_{k-1} - \frac{\V_{k-1} \x_k \x_k^\top \V_{k-1}}{\gamma^{-1} + \x_k^\top \V_{k-1} \x_k}.\label{eq:KalmanUpdateCov}
\end{align}
{Note that at the last iteration, with $k = n$, the Gaussian posterior $p(\thetab|\y_{1:k})$ is perfectly computed with parameters given by Eqs. \eqref{eq:KalmanUpdateMean}-\eqref{eq:KalmanUpdateCov}.} 
The sequence $(\overline{\thetab}_k)_{1 \leq k \leq n}$ turns out to be identical to {the $n$ first iterations of} the {incremental proximal method (IPM)} recursion~\cite{bertsekas2011incremental,bertsekas2011incremental2} applied to Problem \eqref{eqMainProblem}:
\begin{align}\label{eq:IPM}
\overline{\thetab}_k = \prox_{\gamma f_k,\V_{k-1}}(\overline{\thetab}_{k-1}), \quad (\forall k \in \left\{1,\ldots,n\right\})
\end{align}
when {$g=0$ and}
\begin{align}
f_k(\thetab) = \frac{1}{2} (y_k - \x_k^\top \thetab)^2, \quad (\forall \thetab \in \bR^d)
\label{eq:fLS}
\end{align}
for all $k \in \left\{1,\ldots,n\right\}$ and $(\V_k)_{1 \leq k \leq n}$ are specified as in \eqref{eq:KalmanUpdateCov} {(see Props.~4.2--4.4 in \cite{akyildiz2019sequential} for a proof)}
This viewpoint has been extended in \cite{akyildiz2018proximal} for nonlinear least squares, where
\begin{align}
f_k(\thetab) = 
\frac{1}{2}(y_k - h_k(\thetab))^2, \quad \forall \thetab \in \bR^d,
\label{eq:NonLS}
\end{align}
for all $k \in \left\{1,\ldots,n\right\}$. In Eq.~\eqref{eq:NonLS}, each $h_k:\bR^d\to\bR$ is a differentiable function, possibly nonlinear. Thus, the {IPM} iteration of Eq. \eqref{eq:IPM} may not be feasible in a closed form. One can implement the {EKF} for a model with prior
$p(\thetab) = \NPDF(\thetab;\overline{\thetab}_0,\V_0)$ and the likelihood $p(y_k|\thetab) =\NPDF(y_k;h_k(\thetab),\gamma^{-1})$, by linearizing $(h_k)_{1 \leq k \leq n}$. Denoting ${\d}_k = \nabla h_k(\thetab_{k-1})$, we obtain the update rules \cite{akyildiz2018proximal, akyildiz2019sequential}
\begin{align*}
\overline{\thetab}_k &= \overline{\thetab}_{k-1} + \frac{\V_{k-1} {\d}_k (y_k - h_\kappa(\overline{\thetab}_{k-1}))}{\gamma^{-1} + {\d}_k^\top \V_{k-1} {\d}_k}, \\
\V_k &= \V_{k-1} - \frac{\V_{k-1} {\d}_k {\d}_k^\top \V_{k-1}}{\gamma^{-1} + {\d}_k^\top \V_{k-1} {\d}_k}.
\end{align*}
for $k \in \left\{1,\ldots,n\right\}$. {Since the EKF is an approximate Bayesian scheme, multiple passes over the dataset can be performed.}
\section{A Probabilistic IPG method}\label{sec:PIPG}

We now focus on the resolution of the optimization problem in Eq. \eqref{eqMainProblem} when $g \neq 0$. The structure of the cost function suggests the use of the IPG iteration \cite{bertsekas2011incremental,bertsekas2011incremental2}. We consider a variable-metric extension of the IPG. In particular, given \eqref{eqMainProblem}, {the $n$ first iterations of the variable-metric IPG update read as
\begin{align}\label{eq:VarMetIPG}
\overline{\thetab}_k = \prox_{\gamma f_k, \V_{k-1}}(\overline{\thetab}_{k-1} - \gamma \V_{k-1} \nabla g(\overline{\thetab}_{k-1})),
\end{align}
with $\overline{\thetab}_0 \in \mathbb{R}^d$, and $(\V_k)_{k \geq 0} \in \mathbb{R}^{d \times d}$ some predefined SPD matrices. The update \eqref{eq:VarMetIPG} can be viewed as an incremental version of the batch variable-metric PG method that has been extensively studied recently in the optimization literature \cite{Chouzenoux14jota,Combettes_Vu_FB_VM_12}. In the sequel, we propose {a probabilistic interpretation of the IPG which leads to a new update rule for the variable-metric matrices}. We first consider the linear case (i.e., for quadratic $f$ and $g$) for the sake of simplicity, since all computations are tractable {and the inference can be performed in exact manner}. Then we present our general version of the PIPG that encompasses a wider class of cost functions.
\subsection{Linear-Quadratic case}\label{sec:LinQuad}
Let us first assume that $(f_k)_{1 \leq k \leq n}$ is defined {as in} \eqref{eq:fLS} and
\begin{align}
\label{eq:quad}
g(\thetab) = \frac{1}{2} \|\A \thetab\|_2^2 \quad (\forall \thetab \in \mathbb{R}^d),
\end{align}
with $\A \in \mathbb{R}^{m\times d}$, $m \geq 1$. Note that $\A$ is assumed to be known. Using \eqref{eq:quad}, we can write \eqref{eq:VarMetIPG} as
\begin{align}
\widetilde{\thetab}_k &= (\mathbf{I}_d - \gamma \V_{k-1} \A^\top \A) \overline{\thetab}_{k-1},\label{eq:proxGradPred} \\
\overline{\thetab}_k &= \prox_{\gamma f_k,\V_{k-1}}(\widetilde{\thetab}_k),\label{eq:proxGradUpdate}
\end{align}
for $k=1,\ldots,n$. The key observation here is that Eqs.~\eqref{eq:proxGradPred}--\eqref{eq:proxGradUpdate} can be seen as approximate (Kalman) filtering recursions \cite{anderson1979optimal}. To be specific, Eq.~\eqref{eq:proxGradPred} can be seen as the analog to the prediction step within a Kalman filter. Similarly, the update \eqref{eq:proxGradUpdate} can be seen as a Bayesian update using \eqref{eq:fLS}, see Eq.~\eqref{eq:KalmanUpdateMean} \cite{akyildiz2018proximal}. However, Eqs. \eqref{eq:proxGradPred}--\eqref{eq:proxGradUpdate} are different from a Kalman filter, where there would be an update of the covariance matrix between \eqref{eq:proxGradPred}--\eqref{eq:proxGradUpdate}. {Therefore, inspired by Eqs.~\eqref{eq:proxGradPred}--\eqref{eq:proxGradUpdate}, we propose the use of the following state-space model,}
\begin{align}
p(\thetab_0) &= \NPDF(\thetab_0; \overline{\thetab}_0,\V_0), \label{eq:LinModPrior} \\
p(\thetab_k | \thetab_{k-1}) &= \NPDF(\thetab_k; \M_k \thetab_{k-1},\mathbf{0}_{d \times d}), \label{eq:LinModTran} \\
p(y_k|\thetab_k) &= \NPDF(y_k;\x_k^\top \thetab_k, \gamma^{-1}) \label{eq:LinModLik},
\end{align}
where $\mathbf{0}_{d \times d}$, the zero-matrix in $\bR^{d\times d}$, and 
\begin{equation}
\label{eq:M}
\M_k = (\mathbf{I}_d - \gamma \V_{k-1} \A^\top \A) \quad (\forall k \in \left\{1,\ldots,n\right\})
\end{equation}
with $\mathbf{I}_d$ the identity matrix of $\mathbb{R}^d$. 
Now, assume that, the pair $(\overline{\thetab}_{k-1},\V_{k-1})$ is given. We propose to apply filtering recursions for the model \eqref{eq:LinModPrior}--\eqref{eq:LinModLik}, which leads to the PIPG updates for the linear quadratic case. Recursions now consist of a predictive step of the mean and covariance
\begin{align}
\widetilde{\thetab}_k &= \M_k \overline{\thetab}_{k-1},\label{eq:KFMeanPred} \\
\widetilde{\V}_k &= \M_k \V_{k-1} \M_k^\top, \label{eq:KFCovPred}
\end{align}
{respectively, and the update of the mean and covariance},
\begin{align}
\overline{\thetab}_k &= \widetilde{\thetab}_k + \frac{\widetilde{\V}_k \x_k (y_k - \x_k^\top \widetilde{\thetab}_k)}{\gamma^{-1} + \x_k^\top \widetilde{\V}_k \x_k}, \label{eq:KFMeanUpdate}\\
\V_k &= \widetilde{\V}_k - \frac{\widetilde{\V}_k \x_k \x_k^\top \widetilde{\V}_k}{\gamma^{-1} + \x_k^\top \widetilde{\V}_k\x_k} \label{eq:KFCovUpdate},
\end{align}
{respectively,} with $(\M_k)_{1 \leq k \leq n}$ defined in \eqref{eq:M}. It is worth noting that, in Eqs.~\eqref{eq:proxGradPred} and \eqref{eq:proxGradUpdate}, a single $\V_{k-1}$ is used for both iterations. In the corresponding iterations in the proposed method, i.e., {Eqs.~\eqref{eq:KFMeanPred} and \eqref{eq:KFMeanUpdate}}, we make use of $\V_{k-1}$ and $\widetilde{\V}_k$, respectively. In this case, one pass over the dataset is enough since the posterior is exact for the model \eqref{eq:LinModPrior}--\eqref{eq:LinModLik}.
\subsection{General case}
\label{sec:IIIB}
In this section, we present the PIPG algorithm for the general nonlinear case. To be specific, we are going to focus on functions $(f_k)_{1 \leq k \leq n}$ taking the form \eqref{eq:NonLS}. Moreover, we will consider a general function $g$ that we assume to be twice differentiable. In this case, the variable-metric IPG update given in \eqref{eq:VarMetIPG} does not usually yield analytically tractable computations. Moreover, the Kalman recursions, as we presented in the previous section, do not apply. To see this, first consider the mapping $m_\V: \bR^d \mapsto \bR^d$, where
\begin{align*}
m_\V(\overline{\thetab}) = \overline{\thetab} - \gamma \V \nabla g(\overline{\thetab}),
\end{align*}
for some given $\overline{\thetab} \in \bR^d$, $\V \in \bR^{d \times d}$ SPD and $\gamma > 0$. Except when $g$ is quadratic, the above mapping is nonlinear, making it impossible to propagate the uncertainty for the gradient step in \eqref{eq:VarMetIPG} as it was done in \eqref{eq:proxGradPred}. Moreover, when $(f_k)_{1 \leq k \leq n}$ are chosen as in \eqref{eq:NonLS}, it may be complicated to realize the proximal step given in \eqref{eq:VarMetIPG}. To alleviate both problems, we can use the {EKF} \cite{akyildiz2018proximal,anderson1979optimal}. To this end, we build the model
\begin{align}
p(\thetab_0) &= \NPDF(\thetab_0; \overline{\thetab}_0, \V_0),\label{eq:nonlinearPrior} \\
p(\thetab_k | \thetab_{k-1}) &= \NPDF(\thetab_k; m_{\V_{k-1}}(\thetab_{k-1}),\mathbf{0}_{d \times d}),\label{eq:nonlinearTransition}\\
p(y_k | \thetab_k) &= \NPDF(y_k; h_k(\thetab_k),\gamma^{-1}).\label{eq:nonlinearLik}
\end{align}
In order to apply the EKF in the model \eqref{eq:nonlinearPrior}--\eqref{eq:nonlinearLik}, which will lead to \textit{the PIPG algorithm}, we need to linearize the transition model and the observation model. At iteration $k$, given $(\overline{\thetab}_{k-1},\V_{k-1})$ pair, we define the transition matrix,
\begin{align*}
\M_{k} = \mathbf{I}_d - \gamma \V_{k-1} \nabla^2 g(\overline{\thetab}_{k-1}) \quad (\forall k \in \left\{1,\ldots,n\right\}),
\end{align*}
with $\nabla^2 g$ the Hessian map of $g$. Finally, the PIPG updates {can be computed}: first the predicted mean and covariance
\begin{align}
\widetilde{\thetab}_k &= m_{\V_{k-1}} (\overline{\thetab}_{k-1}), \label{eq:EKFMeanPred} \\
\widetilde{\V}_k &= \M_k \V_{k-1} \M_k^\top + \Q, \label{eq:EKFcovPred}
\end{align}
{respectively,\footnote{{Note that, although the dynamical model \eqref{eq:nonlinearTransition} is deterministic {(the process covariance matrix is zero)}, we have introduced $\Q$ in \eqref{eq:EKFcovPred}, a SPD matrix that accounts for the linearization error made by the EKF.}} and then the updated mean and covariance}
\begin{align}
\overline{\thetab}_k &= \widetilde{\thetab}_k + \frac{\widetilde{\V}_k {\d}_k (y_k - h_k(\widetilde{\thetab}_{k}))}{{\gamma}^{-1} + {\d}_k^\top \widetilde{\V}_k {\d}_k}, \label{eq:EKFMeanUpdate}\\
\V_k &= \widetilde{\V}_k - \frac{\widetilde{\V}_k {\d}_k {\d}_k^\top \widetilde{\V}_k}{{\gamma}^{-1} + {\d}_k^\top \widetilde{\V}_k {\d}_k}.\label{eq:EKFCovUpdate}
\end{align}
{respectively, where ${\d}_k = \nabla h_k (\widetilde{\thetab}_k)$. The algorithm is iterated for $k=1,\ldots,n$ and referred to as the PIPG method. 
\begin{rem} Note that Eqs.~\eqref{eq:EKFMeanPred}--\eqref{eq:EKFCovUpdate} are the most general recursions for our method. Like in the linear case presented in Section~\ref{sec:LinQuad}, {sometimes we can simplify} the computations. For instance, if $m_{\V_{k-1}}(\cdot)$ yields a linear mapping for $g$ while $f$ is {a nonlinear least squares loss} as in \eqref{eq:NonLS}, then Eqs.~\eqref{eq:EKFMeanPred}--\eqref{eq:EKFcovPred} simplify into \eqref{eq:KFMeanPred}--\eqref{eq:KFCovPred}. Similarly, when $m_{\V_{k-1}}(\cdot)$ is nonlinear, and $f$ is quadratic as in \eqref{eq:fLS}, then Eqs.~\eqref{eq:EKFMeanUpdate}--\eqref{eq:EKFCovUpdate} simplify into \eqref{eq:KFMeanUpdate}--\eqref{eq:KFCovUpdate}.
\end{rem}
\begin{rem}\label{rem:Differences} 
Although the choice of metrics has been studied in the batch case \cite{Chouzenoux14jota,Becker2011}, no practical ways for choosing them are available in the incremental setting to the best of our knowledge.
The PIPG scheme provides a natural recipe on how to update the metric matrices $(\V_k)_{1 \leq k \leq n}$ in the form of a sequence of posterior covariance matrices.

\end{rem}
\begin{rem} {As mentioned earlier, {in the linear and tractable case the PIPG updates} given by \eqref{eq:KFMeanPred}--\eqref{eq:KFCovUpdate}, are guaranteed to provide, after {$k = n$} iterations { (i.e., after a single pass of the data),} the exact mean and covariance parameters of the Gaussian posterior  associated to the state-space model \eqref{eq:LinModPrior}-\eqref{eq:LinModLik}. {However, the convergence analysis for the general recursions \eqref{eq:EKFMeanPred}--\eqref{eq:EKFCovUpdate} (with inexact Kalman updates) would need further investigation that we leave for future work.} }
\end{rem}

\section{Numerical results}\label{sec:experiments}
In this section, we present two experiments in order to illustrate the performance of PIPG in the context described in Sections~\ref{sec:LinQuad} and \ref{sec:IIIB}.  
{\subsection{Ridge regression}\label{sec:exp:Ridge}
We consider first the linear-quadratic case, depicted in Section~\ref{sec:LinQuad}. We set $\A = \sqrt{\lambda} \mathbf{I}_d$. Moreover, the sought signal $\thetab^\star \in \mathbb{R}^d$ is generated as the realization of a multivariate Gaussian variable using $d = 100$. We then simulated $n = 100,000$ noisy observations $y_k = \x_k^\top \thetab^\star + \eta_k$ with $\eta_k \sim \mathcal{N}(0,1)$ for $k = 1,\ldots,n$.} PIPG recursions \eqref{eq:KFMeanPred}--\eqref{eq:KFCovUpdate} are implemented, for $n = k$ iterations and $40$ step-size values $\gamma$ withing the range $[0.005,0.2]$. We also compare the results to the IPG obtained using \eqref{eq:proxGradPred}-\eqref{eq:proxGradUpdate} with $\mathbf{V}_k = \mathbf{I}_d$ for $k = 1,\ldots,n$. IPG was run with a decaying step-size of the form $\gamma/ k^{0.51}$, for the same range of step-size values than PIPG. Note that running the IPG with a constant step-size causes the algorithm to diverge, therefore we do not show those results. {For both methods, we access the data $(y_k,\x_k)_{1\leq k \leq n}$ in a random order, hence {the} time dependency of $(\x_k)_{1\leq k\leq n}$ is not affecting our results.} We compute the {relative mean squared error} (RMSE) between the current estimate $\thetab_k$ and the true filter coefficient vector $\thetab^\star$ as $\text{E}_k = \|\thetab_k - \thetab^\star\|/\|\thetab^\star\|$. The regularization parameter is set to $\lambda = 10^{-2}$ so as to minimize the final RMSE.}
\begin{figure}[t]
\begin{center}
\hspace*{-0.8cm}
\includegraphics[scale=0.41]{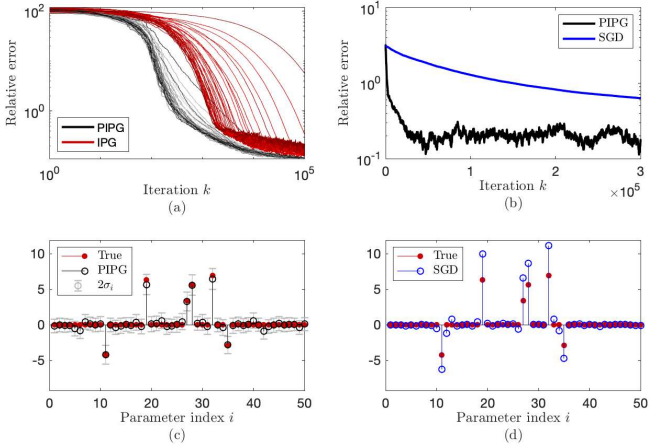}
\end{center}
\caption{(a) Ridge regression example, evolution of RMSE for different runs PIPG and IPG. (b-c-d) Sparse identification example. (b) Evolution of RMSE. (c): Posterior mean and uncertainty estimates, defined as $\pm 2\sigma_i$ for $i = 1,\ldots,d$, $d = 50$, with $(\sigma^2_i)_{1 \leq i \leq 50}$ the diagonal entries of the final posterior covariance matrix $\V_n$, after a single pass over the dataset. (d) SGD estimates after {a single pass} over the dataset.}
\label{Fig:FigureTrueParams}
\end{figure}

{The results are displayed in Fig.~\ref{Fig:FigureTrueParams}(a). It can be seen that PIPG shows {a} stable performance with respect {to} the step-size value. In contrast, IPG appears to be very sensitive to {both} the step-size tuning and {also} the decay rate ({which is} not shown here). Moreover, for a wide range of {step-size} values, PIPG requires less iterations than IPG to achieve minimal RMSE. {Finally, PIPG~provides an estimate of the covariance  as an additional output, which can be particularly useful in practical applications that require an uncertainty quantification in the solution (e.g., biomedical data processing, financial data analytics).}}

\subsection{Sparse nonlinear regression}\label{sec:exp:nonlinear}
Let us now apply the proposed method on the more challenging problem of sparse system identification \cite{diniz2013adaptive,angelosante2010online} under nonlinear observation model. Given a real-valued discrete-time input signal $\big(x_k\big)_{k \in \mathbb{Z}}$, the output of the system at time $k \in \left\{1,\ldots,n\right\}$ is defined as $y_k = h(\x_k^\top \thetab) + w_k$,
where $\x_k = [x_{k-d+1},\ldots,x_k]^\top \in \bR^d$ (assuming circulant boundaries) and $w_k\sim\NPDF(0,\gamma^{-1})$ are i.i.d. measurement noise samples, and $\thetab \in \mathbb{R}^d$ represents the unknown filter taps. A sigmoid nonlinearity $h(u)={1}/(1+\exp(-u))$ for $u\in\bR^d$ is introduced in the system response, modeling for instance some saturation of the sensor. We set the input signal $(x_k)_{k \in \mathbb{Z}}$ as in~\cite{chen2010regularized}, $x_k = a x_{k-1} + \eta_k,$ with $a = 0.8$, $\eta_k \sim \NPDF(0,1)$ and $x_0 \sim \NPDF(0,1)$. We run the PIPG recursions \eqref{eq:EKFMeanPred}--\eqref{eq:EKFCovUpdate} from Section~\ref{sec:IIIB}, where we set $h_k(\thetab) = h(\x_k^\top \thetab)$ for every $\thetab \in \bR^d$, and the regularization function $g$ is chosen as smoothed $\ell_2-\ell_1$ regularization function \cite{chouzenoux2017stochastic}
i.e., $g(\thetab) = \lambda \left(\sum_{i=1}^d \left( 1+ \theta_i^2/\delta^2\right)^{1/2} - 1\right)$ with $\lambda > 0$ and $\delta > 0$ the smoothing parameter. Such regularizer allows to promote sparsity, as when $\delta \to 0$, the $\ell_1$ norm is obtained. {The measurement noise variance is $\gamma^{-1} = 1$}. Note that the parameter $\gamma$ is also the step-size in the proposed method, as we will discuss below. The filter length is $d = 50$ and the output of the system is observed at every time $k\in \{1,\ldots,n\}$ with $n = 300,000$. Regularization parameters are set manually to $(\lambda,\delta) = (10^{-5},0.1)$ so as to reach the best performance in terms of RMSE.  We initialize the PIPG algorithm with a prior distribution with large uncertainty, namely $\V_0 = v_0 \mathbf{I}_d$, where $v_0 = 100$. The process noise covariance matrix, which models the linearization errors in our method, is chosen as $\Q = q \mathbf{I}_d$ with $q = 10^{-4}$. {We set $\gamma = 1$, accordingly with the noise model.} Note that, in general, $\gamma$ is an unknown parameter that is to be set by the user depending on the approximate noise level. For comparison, we implement a stochastic gradient descent (SGD) with learning rate $\gamma^{\textnormal{sgd}}_k = \frac{\alpha_0}{1 + \alpha_1 k}$ for $k = \{1,\ldots,n\}$, where $\alpha_0 = 1$ and $\alpha_1 = 10^{-4}$ which are chosen to reach an optimal decrease. Note that, for this model, it is not possible to implement the IPG since $(f_k)_{1 \leq k \leq n}$ are not easily proximable. 

{Fig.~\ref{Fig:FigureTrueParams}(b) displaying RMSE evolution for both algorithms, shows that the PIPG method {reaches stability in a reduced number of iterations, compared to the SGD, which is a significant practical advantage when one has a limited accessibility to the dataset.}} From Fig.~\ref{Fig:FigureTrueParams}(c)--(d), it can be seen that the PIPG method in Fig.~\ref{Fig:FigureTrueParams}(c) provides a better estimate together with the uncertainty bars $(2\sigma_i)_{1\leq i \leq d}$. A great feature of PIPG is to provide estimates for the covariance matrix, which provides the uncertainty quantification on the parameters. {The behavior of the entries of $(\V_k)_{1\leq k \leq n}$ can be seen from Fig.~\ref{Fig:Covariance}, {along with some comments}}. Let us remark that the computation of this matrix of dimension $d\times d$ implies an increase in computational complexity, as PIPG scales as $\mathcal{O}(d^2)$ while SGD scales as $\mathcal{O}(d)$. 

\begin{figure}[t]
\begin{center}
\includegraphics[scale=0.22]{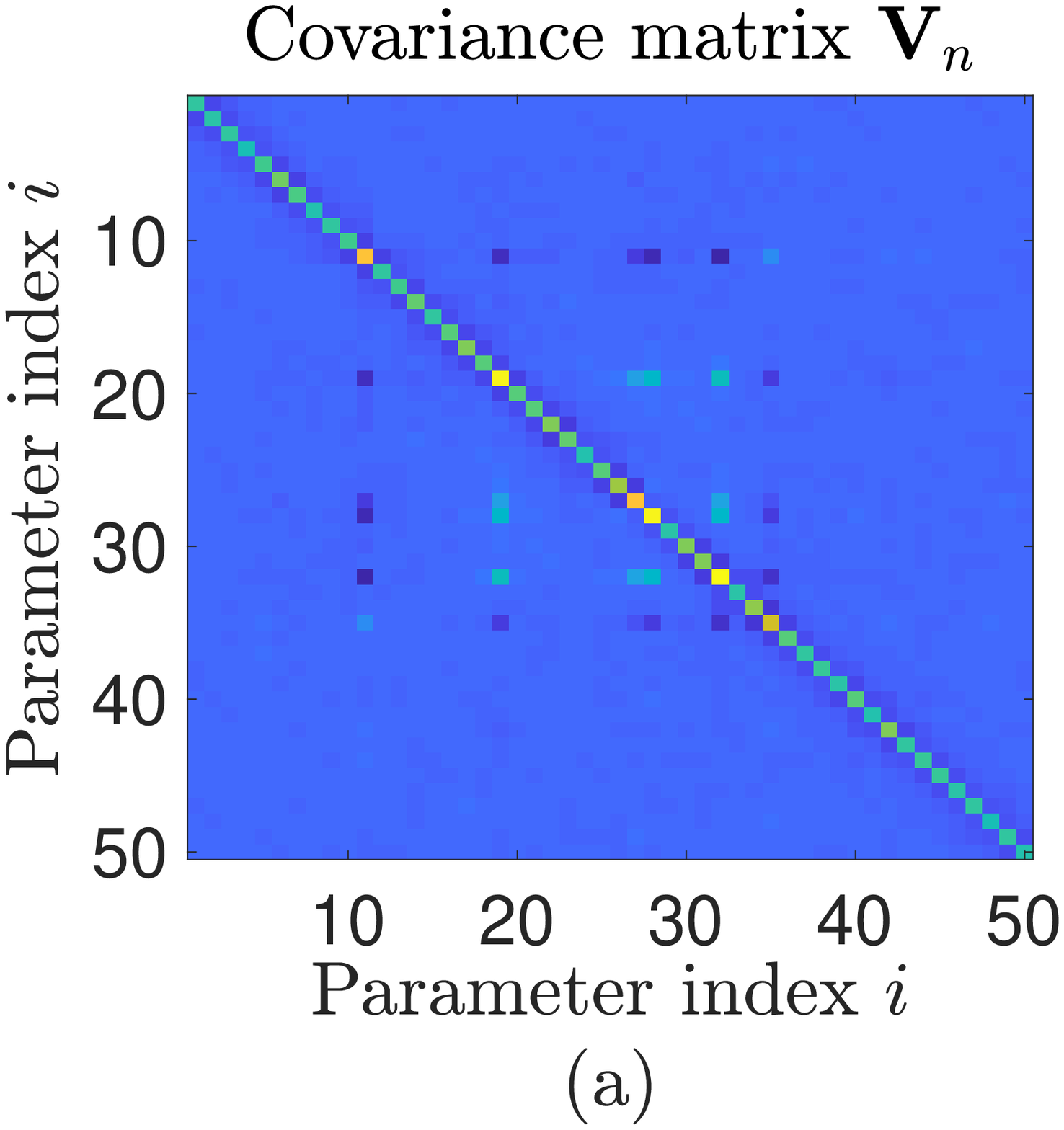}
\includegraphics[scale=0.22]{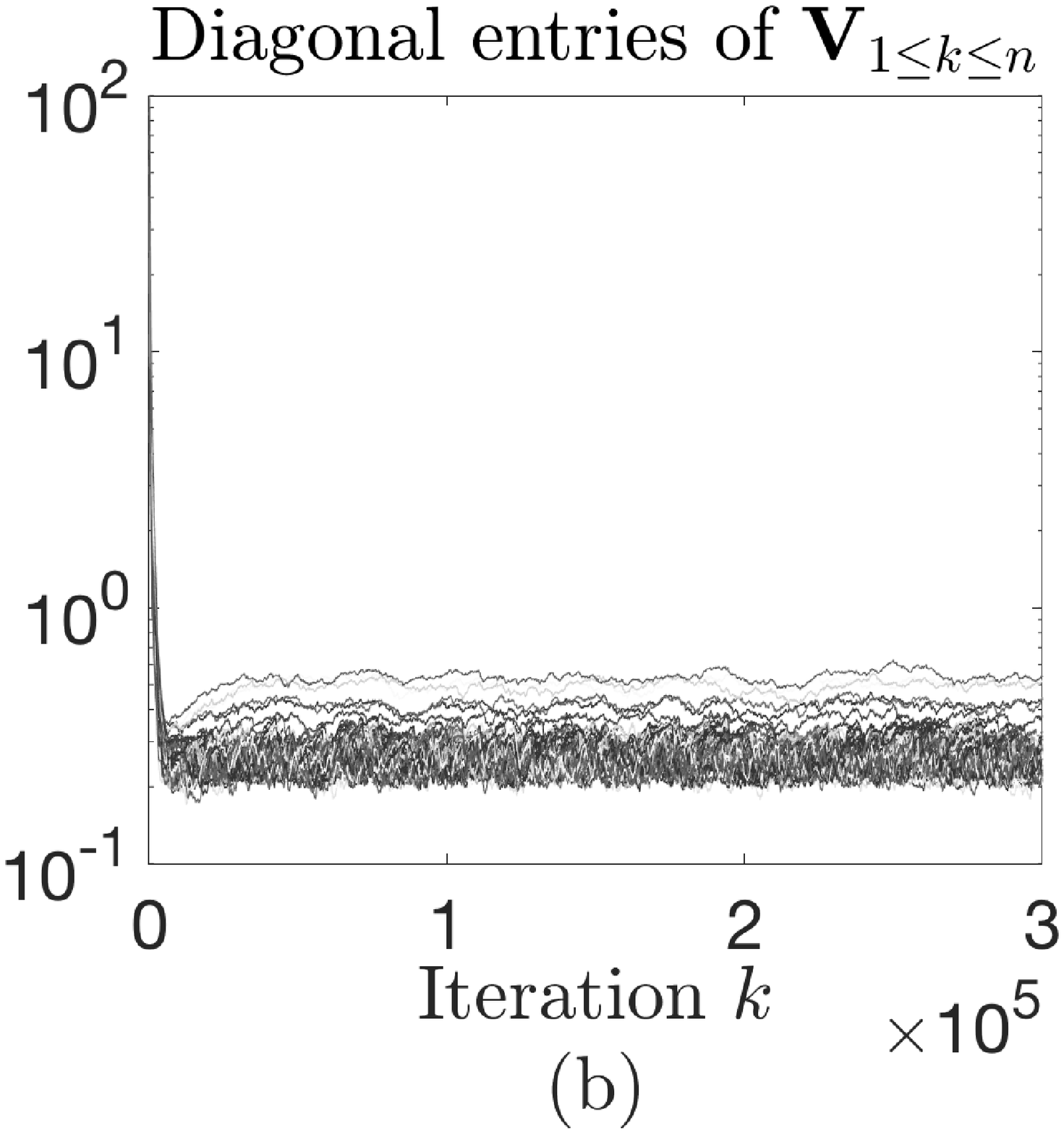}
\end{center}
\caption{The posterior covariance matrix $\V_n$ (a) and the diagonal entries of $(\V_k)_{1 \leq k \leq n}$ (b). From (a), it can be seen that the algorithm learns correlations between dimensions, which enables it to take more efficient steps. From (b), it can be observed that the diagonals of the sequence of covariance matrices converge to certain values quantifying the uncertainty of the final parameters.}
\label{Fig:Covariance}
\end{figure}
\section{Conclusions}\label{sec:conclusions}
We have proposed a probabilistic incremental optimization method which quantifies and propagates the uncertainty over its estimates. In the case of a regularized non-linear least squares, we have reinterpreted the classical IPG method as an approximate inference method in a state-space model. The extension of IPG to the probabilistic setting enables us to provide quantification of the uncertainties inherent in the numerical problem or caused by modeling errors. Our probabilistic interpretation also allows the use of accelerated variable metric updates, whose metric matrices are derived in an automatic and well-defined way. Future investigations will be devoted to {the analysis of the convergence of the PIPG iterates, and the reduction of its complexity by means of suitable approximations}.
\balance
\bibliographystyle{IEEEtran}
\bibliography{draft}

\end{document}